\documentclass[12pt]{article}
\usepackage{theorem}
\usepackage{graphics}
\newtheorem{lemma}{Lemma}

\newtheorem{theorem}[lemma]{Theorem}

{\theorembodyfont{\upshape}}
{\theorembodyfont{\upshape}}
{\theorembodyfont{\upshape}}
{\theorembodyfont{\upshape}}
{\theorembodyfont{\upshape}}
{\theorembodyfont{\upshape}}
\newcommand{\N}{{\bf N}}
\newcommand{\Z}{{\bf Z}}
\newcommand{\R}{{\bf R}}

\newcommand{\rme}{{\rm e}}
\newcommand{\rmd}{\,{\rm d}}

\newcommand{\cH}{{\cal H}}

\newcommand{\cL}{{\cal L}}
\newcommand{\cM}{{\cal M}}
\newcommand{\cP}{{\cal P}}

\newcommand{\sig}{\sigma}
\newcommand{\alp}{\alpha}
\newcommand{\bet}{\beta}
\newcommand{\gam}{\gamma}
\newcommand{\lam}{\lambda}

\newcommand{\eps}{\varepsilon}

\newcommand{\Dom}{{\rm Dom}}

\newcommand{\Spec}{{\rm Spec}}
\newcommand{\Ker}{{\rm Ker}}
\newcommand{\Ran}{{\rm Ran}}
\newcommand{\norm}{\Vert}

\newcommand{\Proof}{\underbar{Proof}{\hskip 0.1in}}
\newcommand{\Note}{\underbar{Note}{\hskip 0.1in}}

\newcommand{\dist}{{\rm dist}}
\newcommand{\lin}{{\rm lin}}
\newcommand{\Schrodinger}{Schr\"odinger }

\newcommand{\pr}{\prime}
\newcommand{\emp}[1]{{\it #1}}

\newcommand{\nopic}[1]{}

\title{SPECTRAL POLLUTION}
\author{E.B. Davies and M. Plum}
\date{9 December 2002}
\begin{document}
\maketitle

\section{Introduction}

It is well known that computing the eigenvalues of a self-adjoint
bounded or differential operator $H$ acting in an infinite
dimensional space $\cH$ is not straightforward. The standard
truncation (or Ritz) method is to take a large finite-dimensional
test function space $\cL$ with an orthonormal basis $\{
\phi_r\}_{r=1}^n$ and to compute the eigenvalues of the $n\times
n$ matrix $A$ given by
\[
A_{r,s}=\langle H\phi_r,\phi_s\rangle.
\]
Unfortunately one often encounters the problem of spectral
pollution: the matrix may possess eigenvalues which are unrelated
to any spectral properties of the original operator. This is an
important problem in several areas of applied mathematics, as is
illustrated by recent papers,
\cite{bg,bdrs,bbg,bdg,boe,densim,mer,rssv,zm}. There are methods
of detecting and resolving such problems, but they are not easy to
understand and for that reason are not well appreciated. Our main
result in this paper is the proof that, at least in some
situations, the geometrically motivated method of Davies,
\cite{dav} is mathematically equivalent to the method of
Zimmermann and Mertins, \cite{zm} which is based on Lehmann's
method \cite{le} and its extension by Goerisch \cite{bg}. We also
provide various other ideas and examples to motivate the method
and show how well it works.

\section{Formulation of the Problem}

We describe the context of the paper at a rather abstract level.
The assumptions are designed to apply to a range of problems,
including those involving second order elliptic operators, as
will become clear in our applications.

Let $A$ be a bounded, self-adjoint, linear operator acting on a
Hilbert space $\cH$, and let $P$ be the orthogonal projection
onto a closed linear subspace $\cL$ of $\cH$. The truncation of
$A$ to $\cL$, denoted $A_{\cL}$, is defined to be the operator
$PAP$ restricted to $\cL$. (Similar definitions may be made for
unbounded $A$ and for quadratic forms.) If $\cL$ is sufficiently
large, one might hope that the spectrum of $A_{\cL}$ is close to
that of $A$. The Rayleigh-Ritz, or variational, theorems
establish precisely this for the part of the spectrum of $A$
outside $[\alp,\bet]$, where $\alp$ (resp. $\bet$) is the smallest
(resp. biggest) point in the essential spectrum of $A$.

Within $[\alp,\bet]$ the situation is far less pleasant. One says
that $\lam$ is a point of spectral pollution of $A$ for the
increasing sequence of subspaces $\cL_n$ with union dense in $\cH$
if there exist $\lam_n\in\Spec(A_{\cL_n})$ such that
$\lam_n\to\lam$ as $n\to\infty$ but $\lam\notin \Spec(A)$. We will
see that spectral pollution is commonplace, but the following two
elementary lemmas sometimes enable one to determine $\Spec(A)$
from its truncations.

\begin{lemma} If $\norm (I-P)AP\norm <\eps$ and
$\lam\in\Spec(A_{\cL})$ then
\[
(\lam-\eps,\lam+\eps)\cap\Spec(A)\not= \emptyset.
\]
\end{lemma}

\Proof One considers $A$ as a perturbation of $PAP+(I-P)A(I-P)$.

\begin{lemma} If $\norm f\norm =1$ and $\norm
Af-\lam f\norm <\eps$ then
\[
(\lam-\eps,\lam+\eps)\cap\Spec(A)\not= \emptyset.
\]
\end{lemma}

\Proof This uses  the resolvent norm estimate
\[
\norm (A-\lam I)^{-1}\norm=\dist \{ \lam,\Spec(A)\}^{-1}.
\]

Lemma 2 is originally due to D. Weinstein. Methods of estimating $\Spec(A)$ starting from a subspace $\cL$
may be divided into two general categories. Prior methods involve
proving that if $\cL_n$ is a particular increasing sequence of
subspaces then the spectrum of $A_{\cL_n}$ converges to that of
$A$. There should also be an explicit error bound, or at least
some information about the rate of convergence. When applied to
differential operators such methods often involve finite element
subspaces and Sobolev space embedding theorems. They have two
disadvantages. The first is that the convergence proofs are often
very difficult and assume conditions which are specific to the
particular applications. The other is that the error bounds are
often far larger than the true errors. Nevertheless valuable and
deep theoretical results can sometimes be obtained using such
methods, \cite{lapsaf1,lapsaf2}.

Posterior methods, on the other hand, involve choosing a subspace
$\cL$ and carrying out computations which provide rigorous bounds
on the location of $\Spec(A)$. If these are good enough one stops,
but if they are not one chooses a different or larger subspace
$\cL$. The problem here is that there is no prior guarantee that
useful bounds will ever be obtained. Nevertheless in practice the
method is often a very good one. Besides the result by D.
Weinstein (Lemma 2), the more refined spectral bounds by T. Kato
\cite{kato}, Aronszajn-A.Weinstein \cite{wei}, Lehmann \cite{le}
and Goerisch \cite{bg} fall into this category. Two particular
posterior methods have recently been the focus of some attention.
One uses the second order spectrum, introduced in \cite{dav,shar},
and is the subject of the companion paper by Levitin and
Shargorodsky, \cite{ls}. This paper studies the other, which
depends upon a systematic study of how certain residuals depend
upon a real parameter, \cite{dav}. We describe a close
relationship between this method and that of Zimmermann and
Mertins, \cite{zm}.

We next set the technical context of the paper. Let $H$ be a
non-negative unbounded self-adjoint operator acting on $\cH$, and
suppose that $0\notin\Spec(H)$. Now suppose that $K$ is a second
non-negative self-adjoint operator, possibly with a different
operator domain, but satisfying a quadratic form relative bound.
In other words if $Q_H$ and $Q_K$ are their two quadratic forms
then $\Dom ( Q_H)\subseteq \Dom(Q_K)$ and
\[
0\leq Q_K (f)\leq aQ_H(f)
\]
for all $f\in \Dom(Q_H)=\Dom(H^{1/2})$. (In fact $Q_K$ need not be
closed, and hence need not be determined by an operator $K$, in
the analysis below.)

The eigenvalue problem $\lam Hf=Kf$ may be recast in the weak, or
variational, form. In other words one seeks $f$, $\lam$ such that
\begin{eqnarray} \label{1} \lam Q_H(f,g)=Q_K(f,g)
\end{eqnarray}
for all $g\in\Dom(Q_H)$.

In many cases the operator $K$ is of the following form. There
exists a closed operator $D:\cH\to\cM$ with operator domain
$\Dom(Q_H)$ such that $H=D^\ast D$ and $K=D^\ast MD$, where $M$
is a non-negative bounded operator acting on $\cM$. Under our
hypotheses $\Ker (D)=0$ and $\Ran (D)$ is closed. There is no loss
of generality in assuming that $\Ran (D)=\cM$, by simply reducing
the size of $\cM$ and restricting $M$. In this case the issue is
to find $f$ such that
\[
\lam \langle Df,Dg\rangle = \langle  M Df,Dg\rangle
\]
for all $g$, or equivalently to find the eigenvalues of the
bounded operator self-adjoint $M$.

Many non-negative differential operators $H$ of order $2n$ can be
written in the form $H=D^\ast D$, where $D$ is a differential
operator of order $n$. Indeed, in a formal sense the assumption is
always satisfied with $D=H^{1/2}$ and $M=H^{-1/2}KH^{-1/2}$. This
observation is rarely useful, because $H^{\pm 1/2}$ are usually
not easy to describe in explicit terms. If $H$ and $K$ are both
differential operators, and $K$ is of lower order than $H$, then
one often finds that $M$ is compact. Its essential spectrum is
then simply $\{0\}$, and the pollution problem does not arise.

Before continuing, we mention that for Sturm-Liouville equations
one can often bypass the problems in this paper by solving the
initial value problem and computing the Prufer angles of the
solution. This can enable one to locate eigenvalues in a gap of
the essential spectrum. The method is not available in higher
dimensions.

\section{Reality of the Phenomenon}

Experience shows that spectral pollution is not just a theoretical
problem. It occurs on a regular basis for operators of physical
importance, such as truncations of periodic \Schrodinger operators
to finite regions; see for example \cite{deihem}. In this section
we demonstrate how bad the phenomenon can be. Later in the paper
we produce more realistic examples.

The following well-known argument (see e.g. \cite{shar}) shows
that spectral pollution can even
occur in an extreme manner for basis functions which are finite
linear combinations of exact eigenfunctions, and which are
therefore in the operator domain.

\begin{lemma} Let $(\lambda_n^{(0)})_{n \in \N}$ and $(\lambda_n^{(1)})_{n
\in \N}$ be two (possibly {\it constant}) sequences of
eigenvalues of problem (\ref{1}),
converging to limits $\lambda^{(0)}$ and $\lambda^{(1)}$,
respectively, satisfying $\lambda^{(0)} < \lambda^{(1)}$.
Moreover, let $\mu$ satisfy $\lambda^{(0)} <\mu <\lambda^{(1)}$.
Then, there exists a sequence $(\psi_n)$ of trial functions which
all are finite linear combinations of exact eigenfunctions of
(\ref{1}), such that
\begin{eqnarray}
\label{2} Q_H (\psi_n, \psi_m) = \delta_{nm}, ~~ Q_K (\psi_n,
\psi_m) = \mu \delta_{nm},
\end{eqnarray}
 i.e., $\mu$ occurs as
the {\it only} eigenvalue of the Ritz method, for every dimension.
\end{lemma}

\Proof Without loss of generality, let $\lambda^{(0)}_n < \mu <
\lambda^{(1)}_n$ for all $n \in \N$, and let $(f_n^{(0)})$ and
$(f_n^{(1)})$ denote sequences of eigenfunctions of (\ref{1})
corresponding to $(\lambda^{(0)}_n)$ and $(\lambda^{(1)}_n)$,
respectively, such that
\begin{eqnarray*}
 Q_H (f_n^{(i)},
f_m^{(i)}) = \delta_{nm}, ~~ Q_K (f_n^{(i)}, f_m^{(i)} ) =
\lambda^{(i)}_n \delta_{nm}.
\end{eqnarray*}
 Finally, choose
$\alpha_n \in [0,2\pi)$ such that
\begin{eqnarray*}
 (\cos^2
\alpha_n) \lambda^{(0)}_n + (\sin^2 \alpha_n) \lambda^{(1)}_n =
\mu, \end{eqnarray*} and \begin{eqnarray*} \psi_n := (\cos
\alpha_n) f_n^{(0)} + (\sin \alpha_n) f_n^{(1)}.
\end{eqnarray*}
 (\ref{2}) is then quite
obvious.

In the same way, also arbitrary {\it sequences} $(\mu_n)$ between
$(\lambda^{(0)}_n)$ and $(\lambda^{(1)}_n)$ can be shown to occur
as spurious eigenvalues.

If we make the weaker assumption that $\lam^{(0)}$ and
$\lam^{(1)}$ are in the essential spectrum of $M$, then a similar
argument can be followed. One can ensure that $D\psi_n$ lie in
the domain of all powers of the operator $M$.

One may criticize the above example on the grounds that the
sequence $(\psi_n)$ of basis functions constructed above may be
far from complete. This gap can be closed as follows, at least if
problem (\ref{1}) has a $Q_H$-orthonormal and complete system of
eigenfunctions:

In this case, let $(f_n^{(2)})_{n \in \N}$ denote a third
sequence of eigenfunctions such that $(f_n^{(0)}), ~
(f_n^{(1)})$, and $(f_n^{(2)})$ together are $Q_H$-orthonormal
and complete. Then, with $(\alpha_n)$ as above, define
\begin{eqnarray*} \varphi_1 := (\cos \alpha_1) f_1^{(0)} + (\sin
\alpha_1) f_1^{(1)}, \end{eqnarray*} and for $n \ge 1$:
\begin{eqnarray*}
\varphi_{4n-2} &:=& (-\sin \alpha_n) f_n^{(0)} + (\cos \alpha_n) f_n^{(1)},\\
\varphi_{4n-1} &:=& (\cos \alpha_{2n}) f_{2n}^{(0)} + (\sin \alpha_{2n}) f_{2n}^{(1)},\\
\varphi_{4n} &:=& (\cos \alpha_{2n+1}) f_{2n+1}^{(0)} + (\sin \alpha_{2n+1}) f_{2n+1}^{(1)},\\
\varphi_{4n+1} &:=& f_{n}^{(2)}. \end{eqnarray*} Then,
$(\varphi_n)_{n \in \N}$ is $Q_H$-orthonormal and complete, and
the Ritz method (with basis functions $\varphi_n$) of dimension
$N$ provides $\mu$ as a $\approx \frac{N}{4}$-fold eigenvalue.

\section{Methods of Avoiding Pollution}

There are various ways of identifying and avoiding spurious
eigenvalues. Suppose that $(\Dom (Q_H),Q_H)$ is a Hilbert space,
with continuous embedding into $\mathcal{H}$. Let $(\cL_n)_{n \in
\N}$ denote an increasing sequence of subspaces, with $\bigcup_{n
\in \N} \cL_n$ dense in $\Dom (Q_H)$. Let $(\mu_n, \chi_n)_{n \in
\N}$ be a sequence of approximate eigenpairs to the problem
(\ref{1}) obtained by the Ritz method:
\begin{eqnarray}\label{3} \chi_n \in \cL_n,~~ \mu_n Q_H (\chi_n,g)
= Q_K (\chi_n,g) \mbox{ for all } g \in \cL_n. \end{eqnarray}
\begin{lemma}
Suppose that $\mu_n \to \mu \in \R$ (or that a subsequence
converges to $\mu$), with $\mu$ not an eigenvalue of (\ref{1});
$\mu$ is therefore either in the essential spectrum or a point of
spectral pollution.

Then, if we normalize $\chi_n$ by \begin{eqnarray}\label{4}
\mbox{\rule[-1mm]{0,5mm}{5mm}} \; \chi_n
\mbox{\rule[-1mm]{0,5mm}{5mm}}_H := \sqrt{ Q_H (\chi_n, \chi_n)} =
1, \end{eqnarray} we have \begin{eqnarray}\label{5} \chi_n
\rightharpoonup 0 \mbox{ weakly in } (\Dom (Q_H),Q_H).
\end{eqnarray}
\end{lemma}

\Note If in particular the embedding $\Dom (Q_H) \hookrightarrow
\mathcal{H}$ is compact, we therefore obtain $\chi_n \to 0$ in
$(\mathcal{H}, \| . \|)$, i.e.,
\begin{eqnarray}\label{6} \frac{\| \chi_n
\|}{\mbox{\rule[-1mm]{0,5mm}{5mm}} \; \chi_n
\mbox{\rule[-1mm]{0,5mm}{5mm}}_H} \to 0. \end{eqnarray} Since in
practice one often performs the Ritz method for various different
dimensions $n$, one has the chance of getting evidence for
(\ref{6}) by evaluating the quotient for various $n$.

\Proof We show that each subsequence of $(\chi_n)$ has a
subsequence converging weakly to zero. Indeed, (\ref{4}) implies
that a subsequence $(\chi_{n_k})$ of a given subsequence converges
weakly in $\Dom (Q_H)$ to some $\chi \in \Dom (Q_H)$. Since $Q_K
\ge 0$ is bounded in $(\Dom (Q_H), Q_H)$, this implies $Q_K
(\chi_{n_k},g) \to Q_K (\chi,g)$ for each $g \in \Dom (Q_H)$.
Consequently, \begin{eqnarray}\label{7} \mu_{n_k} Q_H
(\chi_{n_k}, g) - Q_K (\chi_{n_k},g) \to \mu Q_H (\chi,g) - Q_K
(\chi,g) \end{eqnarray} for each $g \in \Dom (Q_H)$. By
(\ref{3}), the left-hand side of (\ref{7}) is zero for $g \in
\cL_n$ and $n_k \ge n$. Therefore, (\ref{7}) implies
\begin{eqnarray*} \mu Q_H (\chi,g) = Q_K (\chi,g) \mbox{ for all
} g \in \bigcup_{n \in \N} \cL_n \end{eqnarray*} and thus, by
density, for all $g \in \Dom (Q_H)$. Since $\mu$ is no eigenvalue
of (\ref{1}), the assertion $\chi=0$ follows.

One may find the eigenvalues of a bounded operator $M$ in a gap
$(\alp,\bet)$ of its essential spectrum by considering instead
the operator $N=(M-\rho I)^{-1}$, for some $\rho\in(\alp,\bet)$.
The eigenvalues of $M$ in $(\alp,\bet)$ are associated to
eigenvalues of $N$ which lie outside its essential spectrum; the
latter may therefore be found by using the Rayleigh-Ritz
procedure. An apparent disadvantage of this approach is the need
to find an expression for the operator $N$, but this may be
avoided, \cite{zm}.

In this paper we explore a geometrically motivated method of
finding eigenvalues, which traces its origins back to Kato,
\cite{kato}, and Weidmann, \cite[Cor. 6.20]{chat}. One of our main
goals is to describe its close relationship to ideas of Goerisch
and Lehmann, and to developments by various other authors,
\cite{bg, mer, zm}. Suppose that the self-adjoint operator $M$
acting on the Hilbert space $\cM$ has a single eigenvalue $\mu$
which lies within a gap $(\alp,\bet)$ in its spectrum. If
\[
|\lam-\mu|<\min\{ |\lam-\alp|,|\lam-\bet|\}
\]
then one may apply the Rayleigh-Ritz method to the operator
$B=(M-\lam I)^2$ for any $\lam$ close enough to $\mu$. The
smallest eigenvalue of $B$ will then be $(\lam-\mu)^2$, from which
$\mu$ can be calculated (up to a square root ambiguity). The idea
is to choose the value of $\lam$ to maximize the numerical
efficiency of the method. One's intuition is to take $\lam$ as
close as possible to the eigenvalue $\mu$, but we shall see that
this is not the correct strategy.

As in \cite{dav}, we investigate the $\lam$-dependence of the
procedure above by introducing the function
\[
F(\lam)=\min\{\norm Mf-\lam f\norm:f\in\cM\mbox{ and }\norm
f\norm =1\} .
\]
If one could evaluate $F$ numerically then the formula
\begin{equation}
F(\lam)=\dist\{\lam,\Spec(M)\}\label{Fdist}
\end{equation}
would enable one to locate the spectrum of $M$. To find the
eigenvalues of $M$ one needs to combine this with some
theoretical input which distinguishes between the point spectrum
and essential spectrum of $M$.

Following \cite{dav}, we introduce the approximate functions
\[
F_n(\lam)=\min\{\norm Mf-\lam f\norm:f\in \cL_n\mbox{ and }\norm
f\norm =1\}.
\]
These may be used to obtain useful spectral information, as
follows.

\begin{lemma}
If $M$ is bounded then for any $\eps >0$ there exists $N_\eps$
such that $n\geq N_\eps$ implies
\begin{equation}
F(s)\leq F_n(s)\leq F(s)+\eps\label{conv}
\end{equation}
for all $s\in\R$.
\end{lemma}

\Proof It follows from their definition that $F_n(s)$ decrease
monotonically to $F(s)$ as $n\to\infty$ for each $s\in\R$. It
also follows from its definition that
\[
|F_n(\lam)-F_n(\mu)|\leq |\lam-\mu|
\]
for all $\lam,\mu\in\R$. Equivalently $|F^\pr(\lam)|\leq 1$ for
almost all $\lam\in\R$ (in the sense of Lebesgue). Since the
functions $F_n$ are equicontinuous they must converge uniformly
on compact subsets of $\R$. Suppose that $n\geq N_\eps$ implies
(\ref{conv}) for $-\norm M\norm \leq s \leq \norm M\norm$. If
$s>\norm M \norm$ then
\begin{eqnarray*}
F(s)&\leq &F_n(s)\\
&\leq & F_n(\norm M \norm )+s-\norm M \norm\\
&\leq& F(\norm M \norm)+\eps +s-\norm M \norm\\
&=& F(s)+\eps.
\end{eqnarray*}
A similar argument applies if $s< -\norm M \norm$.

\begin{lemma}
Given (\ref{conv}), let $\lam$ be an eigenvalue of $M$ and let
there be no other point of $\Spec(M)$ in
$[\lam-4\eps,\lam+4\eps]$. Then $F_n$ has a local minimum at some
point  $\sig\in [\lam-\eps,\lam+\eps]$. If $\mu=\sig-F_n(\sig)$
and $\nu=\sig+F_n(\sig)$ then
\[
\lam-2\eps\leq \mu \leq \lam\leq \nu\leq\lam+2\eps.
\]
\end{lemma}

\Proof Since $F(\lam)=0$ we see that $0\leq F_n(\lam)\leq \eps$.
However if $\eps<|s-\lam|\leq 2\eps$ we have
\[
\eps<|s-\lam|=F(s)\leq F_n(s).%
\]
Hence the minimum of $F_n$ within the interval
$[\lam-2\eps,\lam+2\eps]$  must occur for $\sig$ within
$[\lam-\eps,\lam+\eps]$. Moreover $0\leq F_n(\sig)\leq \eps$. The
second statement of the lemma follows immediately.

The functions $F_n(\cdot)$ may be computed as follows. One
introduces
\[
F_{m,n}(\lam)=\min\{\norm P_m(Mf-\lam f)\norm:f\in \cL_n\mbox{ and
}\norm f\norm =1\}
\]
where $m\geq n$. If $m=n$ then this achieves no more than finding
the eigenvalues of the truncation of $M$ to $\cL_n$ (i.e. applying the Ritz method), but if $m\gg
n$ it may be far superior. One may also compute $F_n(\lam)$
exactly by replacing $\cL_m$ above by the finite-dimensional
subspace
\[
\tilde \cL_n=M(\cL_n)+\cL_n.
\]
The projection $P_m$ in the argument below is then replaced by
the projection $\tilde P_n$ onto $\tilde \cL_n$. Whether or not this
Lanczos process is easily implemented depends upon the example.

The evaluation of $F_{m,n}(\lam)$ is straightforward. One has the
formula
\[
F_{m,n}(\lam)=G(m,n,\lam)^{1/2}
\]
where $G(m,n,\lam)$ is the smallest eigenvalue of the $n\times n$
matrix
\[
B(m,n,\lam)=P_n M P_m M P_n-2\lam P_n M P_n +\lam^2 P_n.
\]

In practice one finds potential eigenvalues $\lam$ of $M$ by using
the Ritz method. One then examines $F_{m,n}(\lam)$ for those
values. If it is sufficiently small then one has good evidence
that there is indeed an eigenvalue near the candidate value,
while if it is not small then one has no such evidence. Using
only a finite number of matrix entries of the operator $M$ no
more could ever be established, but this method is far better
than simply using the Ritz method without any checks.

The above provides evidence for the existence of an eigenvalue
near a specified number. One can also obtain evidence for its
multiplicity, or more precisely for the number of eigenvalues in a
very short interval, by an elaboration of the above ideas which is
described in \cite{dav}.

\section{Bounding the eigenvalues}

In \cite{dav} we discussed how to obtain rigorous bounds on the
spectrum of $M$ (or of problem (\ref{1}), respectively) from the functions $F_n(\cdot)$, thus avoiding the
problem of spectral pollution. In this paper we develop the ideas
further. If the subspace $\cL_n$ is large enough then $F_n(\lam)$
has local minima near the eigenvalues of $H$, and its derivatives
are close to $\pm 1$ on either side of the eigenvalues. The
following exact formula for the derivative may be useful, for
example if applying Newton's method. We have
\[
F_n^\pr(\lam)=\frac{G_n^{\,\pr}(\lam)}{2\sqrt{G_n(\lam)}}
\]
where $G_n(\lam)$ is the smallest eigenvalue of the $n\times n$
matrix
\begin{eqnarray}
B_n(\lam)&=&P_n M^\ast M P_n-2\lam P_n M P_n +\lam^2 P_n\\
&=& D_n-2\lam M_n+\lam^2I\label{B}
\end{eqnarray}
where the second line is understood as acting in $\cL_n$ and
$M_n=M_{\cL_n}$. If $f_n(\lam)$ is the normalized eigenvector
corresponding to the eigenvalue $G_n(\lam)$ then it is well-known
that
\begin{eqnarray}
G_n^{\,\pr}(\lam)&=&\langle B_n^\pr(\lam)f_n(\lam),f_n(\lam)\rangle\\
&=&2\lam -2\langle M_n f_n(\lam),f_n(\lam)\rangle.
\end{eqnarray}

A computation of the graph of the function $F_n(\cdot)$ is bound
to be slow if $n$ is large. It is best therefore to attempt this
only for reasonably small $n$, in order to understand
approximately where the spectrum is located. One then passes to
larger values of $n$ for which $F_n(\lam)$ is only computed for
carefully selected values of $\lam$.

We next discuss the computation of the graph of $F_n(\cdot)$ for
reasonably small $n$. The procedure below is guided by the fact
that it is optimal if $D_n=M_n^2$, or equivalently if $(I-P_n) M
P_n=0$. It should therefore be close to optimal if
$D_n-M_n^2=P_nM(I-P_n)MP_n$ is small. The starting point is to
find all of the eigenvalues $\lam_r$ of $M_n$ and then to compute
$G(\lam_r)$ for each $r$. If $G(\lam_r)$ is not sufficiently small
for some $r$ then $\lam_r$ is rejected as spurious. Let the
remaining eigenvalues be listed in increasing order. The
corresponding normalized eigenvectors $g_r$ of $M_n$ are used as
starting vectors for an inverse power iteration to find
$f_r=f(\lam_r)$. As a check it should be verified that $f_r$ and
$g_r$ are approximately equal, and that the set $f_r$ is
approximately orthonormal as $r$ varies.

Now assume that one wants a list of values of $F_n(\cdot)$ over a
set of points covering the interval $(\lam_r,\lam_{r+1})$, in
order to plot a graph. This is done separately on either side of
the midpoint $\gam_r$ of the interval, and we consider the left
hand half. One chooses some subdivision
$\lam_r=\sig_1<\sig_2<...<\sig_m<\gam_r$ and then computes
$G(\sig_r)$ in order. At each stage the starting vector for the
inverse iteration should be the eigenvector for the previous
stage. In the right hand half interval one deals with the points
in decreasing order.

The reason for using different procedures on either side of the
midpoint becomes clear if one considers the exact case, in which
$D_n=M_n^2$; the eigenvector is then constant in each subinterval
and changes discontinuously as $\lam$ increases through $\gam_r$.
In the inexact case one should expect $f(\lam)$ to change rapidly
near $\lam=\gam_r$, and it is therefore advisable to start the
iteration with the vector $f(\lam^\pr)+f(\lam^{\pr\pr})$, where
$\lam^\pr$ and $\lam^{\pr\pr}$ are values already computed on
either side of $\gam_r$.

One also needs to evaluate $F_n(\lam)$ for large values of $n$ at
individual points $\lam$ satisfying $\lam_r <\lam<\lam_{r+1}$.
Once again this may be done using inverse power iteration, but
using $f_r+f_{r+1}$ as the initial vector.

We finally discuss how the computed values of $F_n(\lam)$ are
used to locate the spectrum of $M$. We use
\[
\Spec(M)\cap [\lam-F_n(\lam),\lam+F_n(\lam)]\not=\emptyset
\]
for all $\lam\in\R$, which follows from the spectral theorem, and
the obvious
\[
F_n(\lam)\geq F(\lam)
\]
for all $\lam\in\R$. Unless $F_n(\lam)$ and $F(\lam)$ are
approximately equal, the above provide little information. One
has two choices:

(i) One can supplement the computation of $F_n(\lam)$ by prior
information about the location of the spectrum of $M$.

(ii) One can examine the graph of $F_n(\lam)$ to see whether its
form is consistent with the assumption that it is approximately
equal to $F(\lam)$, using (\ref{Fdist}). If not then one can
increase the value of $n$, or revert to (i).

Of course the former procedure is to be preferred, but the latter
is likely to be easier. We illustrate both in the example
discussed in the next section.

For the rest of the paper we make the following hypothesis.

{\it (S) We have $F_n(s)>F(s)$ and $|F^\pr_n(s)|<1$ for all
$s\in\R$.}

This can be ensured by replacing $F_n$ by
\[
\hat F_n(s)=F_n(s)+\eps\rme^{-F_n(s)}
\]
for any $\eps>0$. Increasing $F_n(s)$ slightly makes the spectral
bounds which we prove slightly weaker, but does not involve an
essential loss of generality. We will use the following
consequence of the hypothesis (S).
\[
\Spec(M)\cap (\lam-F_n(\lam),\lam+F_n(\lam))\not=\emptyset.
\]

Our goal is to determine those eigenvalues of $M$ which lie in a
given interval $(\alp,\bet)$; $\alp$ and $\bet$ might sometimes be
known points of the essential spectrum of $M$. Let $\sig_r$,
$1\leq r\leq R$ be an increasing sequence of numbers lying in the
interval $(\alp,\bet)$; they are our first approximations to the
eigenvalues of $M$ within the interval. They might be the
non-spurious eigenvalues $\lam_r$ of $M_n$ already described, or
they might be the local minima of the function $F_n(\lam)$,
obtained by solving $F_n^\pr(\lam)=0$ using the bisection method.
In the simplest cases these two methods of choosing $\sig_r$ yield
approximately the same values. If this fails one can still
proceed, using the second choice of $\sig_r$ just mentioned.
Putting $\mu_r=\sig_r - F_n(\sig_r)$ and
$\nu_r=\sig_r+F_n(\sig_r)$ for $1\leq r\leq R$, we have
\[
\Spec(M)\cap (\mu_r,\nu_r)\not=\emptyset
\]
for such $r$. We also put $\nu_0=\alp$ and $\mu_{R+1}=\bet$.
Further progress depends upon the following hypothesis. Following
(i) above, this might have been proved, possibly using a homotopy
method. Or following (ii), it might simply be motivated by the
numerical evidence, and tested by the consistency of the
consequences drawn using it.

{\it(H) We have
\[
\nu_0<\mu_1<\nu_1<\mu_2<...<\mu_R<\nu_R<\mu_{R+1}.
\]
The intervals $(\nu_{r-1},\mu_r)$ do not meet $\Spec(M)$ for
$1\leq r\leq R+1$, while $[\mu_r,\nu_r]$ each contain a single
eigenvalue $m_r$ of $M$ for $1\leq r\leq R$.}

Situations in which the operator $M$ has eigenvalues which are too
close to be resolved using the subspace $\cL_n$ can often be
handled using a modified procedure, or more simply by increasing
$n$.

We describe an iterative procedure for obtaining steadily more
accurate enclosures of the eigenvalues. The following lemma
provides the starting point for the iteration. The constructions
are most easily understood by plotting the graph of a typical
function $F_n$. See Figure 1 for an example in which $\alp=0$,
$\bet=1$ are eigenvalues of infinite multiplicity, and $R=1$.

\begin{lemma}\label{induct} Assuming (S) and (H), we have
\[
\mu_r < \mu_r^{\pr}< m_r < \nu_r^{\pr} < \nu_r
\]
for $1 \leq r\leq R$, where $\nu_r^{\pr}$ and $\mu_r^{\pr}$ are
defined by solving the following equations.
\begin{eqnarray*}
F_n(s_r)&=& s_r-\nu_{r-1}\\
\nu_r^{\pr}&=&s_r+F_n(s_r)\\
F_n(t_r)&=&\mu_{r+1}-t_r\\
\mu_r^{\pr}&=&t_r-F_n(t_r)
\end{eqnarray*}
\end{lemma}

\Proof Let $1\leq r\leq R$. The function
\[
\phi(s)=F_n(s)-s+\nu_{r-1}
\]
satisfies $\phi^\pr(s)=F_n^\pr(s)-1< 0$ for all $s$,
$\phi(\nu_{r-1})=F_n(\nu_{r-1})>0$ and
\[
\phi(\sig_r)=F_n(\sig_r)-\sig_r+\nu_{r-1}=\nu_{r-1}-\mu_r <0.
\]
Hence $\nu_{r-1}<s_r<\sig_r$. We next observe that
\[
F_n(s_r)-F_n(\sig_r)< \sig_r-s_{r}
\]
may be rewritten in the form $ \nu_{r}^{\pr}< \nu_r$. Also
\[
\Spec(M)\cap (s_{r}-F_n(s_{r}),s_{r}+F_n(s_{r}))\not=\emptyset
\]
may be rewritten as
\[
\Spec(M)\cap (\nu_{r-1},\nu_{r}^{\pr})\not=\emptyset.
\]
Using (H) we now deduce that $m_r < \nu_{r}^{\pr}$. The proof
that $\mu_r < \mu_r^{\pr} < m_r$ is similar.

Improving the above result iteratively depends upon setting up
the inductive hypothesis properly. Our hypothesis below allows the
possibility of increasing the value of $n$ from one stage of the
iteration to the next, in other words of increasing the size of
the test function space as the iteration proceeds. In practice one
would choose a fairly small initial value of $n$ and apply the
lemma repeatedly until the iteration stops yielding improvements,
and then increase $n$ significantly and repeat the iteration with
that value of $n$ until improvements cease.

\begin{theorem}
Let $\mu_r$, $\mu_r^\pr$, $\nu_r$, $\nu_r^\pr$, $s_r$ and $t_r$
satisfy
\[
\mu_r\leq \mu_r^{\pr}< m_r < \nu_r^{\pr}\leq \nu_r
\]
and
\begin{eqnarray*}
F_n(s_r)&\leq& s_r-\nu_{r-1}\\
\nu_r^{\pr}&\geq&s_r+F_n(s_r)\\
F_n(t_r)&\leq&\mu_{r+1}-t_r\\
\mu_r^{\pr}&\leq&t_r-F_n(t_r)
\end{eqnarray*}
for $1\leq r\leq R$. Put $\nu_0^\pr=\sig_0$ and
$\mu_{R+1}^\pr=\sig_{R+1}$. Given (S) and (H), define
$\mu_r^{\pr\pr}$, $\nu_r^{\pr\pr}$, $s_r^\pr$ and $t_r^\pr$ for
$1\leq r\leq R$ by
\begin{eqnarray*}
F_n(s_r^\pr)&=& s_r^\pr-\nu_{r-1}^\pr\\
\nu_r^{\pr\pr}&=&s_r^\pr+F_n(s_r^\pr)\\
F_n(t_r^\pr)&=&\mu_{r+1}^\pr-t_r^\pr\\
\mu_r^{\pr\pr}&=&t_r^\pr-F_n(t_r^\pr)
\end{eqnarray*}
Then
\[
\mu_r\leq \mu_r^{\pr}\leq\mu_r^{\pr\pr} < m_r <
\nu_r^{\pr\pr}\leq\nu_r^{\pr}\leq \nu_r
\]
for $1\leq r\leq R$.
\end{theorem}

\Proof Let $1\leq r\leq R$. The function
\[
\phi(s)=F_n(s)-s+ \nu_{r-1}^\pr
\]
satisfies $\phi^\pr(s)<0$ for all $s\in\R$ and
\[
\phi(\nu_{r-1}^\pr)=F_n(\nu_{r-1}^\pr)>0.
\]
Also
\[
\phi(s_r)=F_n(s_r)-s_r+\nu_{r-1}^\pr\leq \nu_{r-1}^\pr -\nu_{r-1}
\leq 0.
\]
Hence $\nu_{r-1}^\pr <s_r^\pr \leq s_r$. We next observe that
\[
F_n(s_r^\pr)-F_n(s_r)\leq s_r-s_r^\pr
\]
implies
\[
\nu_r^{\pr\pr}=F_n(s_r^\pr)+s_r^\pr\leq F_n(s_r)+s_r\leq
\nu_r^\pr .
\]
Also
\[
\Spec(M)\cap
(s_r^\pr-F_n(s_r^\pr),s_r^\pr+F_n(s_r^\pr))\not=\emptyset
\]
is equivalent to
\[
\Spec(M)\cap (\nu_{r-1}^\pr,\nu_r^{\pr\pr})\not=\emptyset
\]
and implies that $m_r<\nu_r^{\pr\pr}$. The proof that
$\mu_r^\pr\leq \mu_r^{\pr\pr}< m_r$ is similar.

At each stage in the implementation of the iterative process one
has further consistency checks, namely the numerical procedure
need not yield the inequality $\mu_r^{\pr\pr}< \nu_r^{\pr\pr}$
for all $r$. If this does not happen then (H) must be false; that
is the operator $M$ must have more spectrum than assumed.

Plotting the graph of a typical function $F_n$ suggests certain
results about the location of $s_r$ and $t_r$; these are proved
below.

\begin{lemma}\label{sconv} Assume (S) and (\ref{conv}).
Given $\lam\in\R$, let $\lam^\pr$ denote the smallest point in
$\Spec(M)$ greater than $\lam$.  Assume $\lam^\pr>\lam+4\eps$ and
let $s$ be the solution of $F_n(s)=s-\lam$. Then
\[
\frac{\lam+\lam^\pr}{2} < s < \frac{\lam+\lam^\pr}{2}+2\eps.
\]
\end{lemma}

\Proof We combine the fact that the function
$\phi(s)=F_n(s)-s+\lam$ is strictly monotone decreasing and
continuous with the inequalities
\[
\phi((\lam+\lam^\pr)/2)>F((\lam+\lam^\pr)/2)-(\lam+\lam^\pr)/2+\lam=0\\
\]
and
\[
\phi((\lam+\lam^\pr)/2+2\eps)<F((\lam+\lam^\pr)/2+2\eps)
+2\eps-(\lam+\lam^\pr)/2-2\eps+\lam=-2\eps< 0.
\]

\begin{theorem} Suppose that $M$ has a single eigenvalue $m$ in $(\alp,\bet)$.
Define $s_n$ and $t_n$ by
\begin{eqnarray*}
F_n(s_n)&=& s_n-\alp\\
\nu_n&=& s_n+F_n(s_n){\hskip 0.25 in} (=2s_n-\alp)\\
F_n(t_n)&=& \bet-t_n\\
\mu_n&=& t_n-F_n(t_n){\hskip 0.25 in} (=2t_n-\bet).
\end{eqnarray*}
Then
\[
\mu_n<m<\nu_n
\]
for all $n\geq 1$. Moreover
\begin{equation}
\lim_{n\to\infty}s_n= (\alp+m)/2 ,{\hskip 0.2in}
\lim_{n\to\infty}t_n=(m+\bet)/2\label{convrate}
\end{equation}
and
\[
\lim_{n\to\infty}\mu_n=\lim_{n\to\infty}\nu_n =m.
\]
\end{theorem}

\Proof The first statement follows from Lemma \ref{induct} with
$R=1$, $\nu_0=\alp$ and $\mu_2=\bet$. The other statements are
then consequences of Lemma \ref{sconv}.

\section{The Zimmermann-Mertins Paper}

We have repeatedly referred to solving the equations
$F_n(s)=s-\nu$ and $F_n(t)=\mu-t$ for various $\mu,\nu$. This may
be carried out by bisection or by Newton's method. By examining
the paper of Zimmermann and Mertins, \cite{zm}, one discovers a
method which is probably better in most circumstances; see
(\ref{K}) below.

We assume that $(\nu,\mu)$ is an interval which contains a single
eigenvalue $m$ of $M$, but no other point of Spec $(M)$. E.g. in the situation described in hypothesis (H), we have in mind tho choose $\nu = \nu_{r-1}, \mu = \mu_{r+1}, m = m_r$, for any $r \in \{ 1, \dots, R\}$. Suppose
for simplicity that $\nu$ and $\mu$ are not eigenvalues of
$M$.\bigskip

Let $\cL_n \subset \mathcal{M}$ denote an $n$-dimensional trial
function space. To ensure that the exact eigenfunction
corresponding to $m$ is `not too far' from $\cL_n$, we make
the additional assumption
\begin{eqnarray}
\nu < \max\limits_{v \in \cL_n \setminus \{ 0 \}}  \frac{\langle
Mv,v \rangle}{\langle v,v \rangle},~~ \min\limits_{v \in \cL_n
\setminus \{ 0 \} }   \frac{\langle Mv,v \rangle}{\langle v,v
\rangle} < \mu.\label{A}
\end{eqnarray}
(Without (\ref{A}), the Lehmann method gives nothing, and the
$F_n$-method gives a lower bound for $m$ less than $\nu$, resp.
an upper bound larger than $\mu$, which makes it more or less
useless.)\bigskip

The right-definite Lehmann method may be applied to the situation
in which only one eigenvalue is to be bounded. We follow
\cite[Theorem 4.2]{zm}, where the spectral parameter $\rho$ is
once chosen to be $\mu$, and once $\nu$. Let $\{ v_1, \dots,
v_n\}$ denote a basis of $\cL_n$, and let $\tau^-$ be the smallest
eigenvalue of the $n \times n$ matrix eigenvalue problem
\begin{eqnarray}
( \langle M v_i - \mu v_i,v_k \rangle )_{i,k} \;x = \tau (\langle
M v_i - \mu v_i,Mv_k - \mu v_k  \rangle )_{i,k}\; x \label{B2}
\end{eqnarray}
and $\tau^+$ the largest eigenvalue of
\begin{eqnarray}
( \langle M v_i - \nu v_i,v_k \rangle )_{i,k} \;x = \tau (\langle
M v_i - \nu v_i,Mv_k - \nu v_k  \rangle )_{i,k} \;x. \label{C}
\end{eqnarray}
Then
\begin{eqnarray}
\mu  + \frac{1}{\tau^-} \le m \le \nu + \frac{1}{\tau^+}.
\label{D}
\end{eqnarray}
(Note that $\tau^- < 0 < \tau^+$, according to (\ref{A})).

The following theorem states that the $F_n$-method gives exactly
the right-definite Temple-Lehmann-Maehly bounds, in the
Zimmermann-Mertins setting. See also \cite[Theorem 6]{dav}.

\begin{theorem} Let $\tau^-$ and $\tau^+$ as above, and let $s$ and
$t$ denote the solutions of (cf. Lemma \ref{induct}):
\begin{eqnarray}
F_n (s) = s - \nu, ~~~F_n (t) = \mu - t.\label{E}
\end{eqnarray}
Then,
\begin{eqnarray}
s + F_n (s) = \nu + \frac{1}{\tau^+}, ~~~t-F_n (t) = \mu +
\frac{1}{\tau^-}. \label{F}
\end{eqnarray}
Hence
\begin{eqnarray}
s=\nu+\frac{1}{2\tau^+}, ~~~t=\mu+\frac{1}{2\tau^-}.\label{K}
\end{eqnarray}
\end{theorem}

\Proof For $\rho \in \R$, define
\begin{eqnarray*}
A_2 (\rho) &:=& ( \langle M v_i - \rho v_i,M v_k - \rho v_k \rangle )_{i,k=1,\dots,n}, \\
A_1 (\rho) &:=& ( \langle M v_i - \rho v_i, v_k \rangle
)_{i,k=1,\dots,n}, ~~A_0 := (\langle v_i, v_k \rangle )_{i,k = 1,
\dots, n}.
\end{eqnarray*}
Straightforward calculations and use of (\ref{E}) give, for all $x
\in \R^n$,
\begin{eqnarray}
x^T A_2 (s) x - F_n (s)^2 x^T A_0 x &=& x^T A_2 (\nu) x + 2 (\nu -s) x^T A_1 (\nu) x \nonumber \\
&&~~+ [(\nu-s)^2 - F_n (s)^2] x^T A_0 x  \nonumber \\
&=& x^T A_2 (\nu) x - 2 F_n (s) x^T A_1 (\nu) x.\label{G}
\end{eqnarray}
By the definition of $F_n (s)$, the left-hand side of (\ref{G}) is
non-negative for all $x$, and zero for some specific $x$. This is
therefore true also for the right-hand side, i.e.,
\begin{eqnarray*} \frac{x^T A_1 (\nu) x}{x^T A_2 (\nu)x} \le
\frac{1}{2 F_n (s)} ~~{\rm ~for~all~}x \in \R^n \setminus \{ 0
\}, \end{eqnarray*} with equality for some specific $x$. Thus,
$1/(2F_n (s))$ is the maximum of the Rayleigh quotient of problem
(\ref{C}), i.e.,
\begin{eqnarray}
\frac{1}{2 F_n (s)} = \tau^+.\label{H}
\end{eqnarray}
Using (\ref{E}) again, we obtain
\begin{eqnarray*}
s + F_n (s) = \nu
+ 2 F_n (s) = \nu + \frac{1}{\tau^+}.
\end{eqnarray*}

A similar method yields
\begin{eqnarray*}
x^T A_2 (t) x - F_n (t)^2 x^T A_0 x = x^T A_2 (\mu) x + 2 F_n (t)
x^T A_1 (\mu) x,
\end{eqnarray*}
and the same arguments as before yield
\begin{eqnarray*}
\frac{x^T A_1 (\mu) x}{x^T A_2 (\mu)x} \ge -\frac{1}{2 F_n (t)}
~~{\rm ~for~all~}x \in \R^n \setminus \{ 0 \},
\end{eqnarray*}
with equality for some specific $x$, which implies
\begin{eqnarray}
- \frac{1}{2 F_n (t)} = \tau^- \label{J}
\end{eqnarray}
and thus
\begin{eqnarray*}
t - F_n (t) = \mu - 2 F_n (t) = \mu +
\frac{1}{\tau^-}.
\end{eqnarray*}

The final statement of the theorem is obtaining by combining
(\ref{E}), (\ref{H}) and (\ref{J}).

\section{An Example in One Dimension}

Suppose that the operator $H$ has a complete orthonormal set of
eigenfunctions $\{ \phi_r\}_{r=1}^\infty$, whose corresponding
eigenvalues $\lam_r$ are all positive. Suppose also that the
eigenvalues are listed in increasing order and tend to infinity.
If we put
\begin{equation}
\psi_r=\lam_r^{-1/2}D\phi_r   \label{basis}
\end{equation}
then $\{\psi_r\}_{r=1}^\infty$ forms an orthonormal basis in
$\cM$ and the matrix elements of $M$ with respect to this basis
are
\[
M_{r,s}=\lam_r^{-1/2}\lam_s^{-1/2}Q_K(\phi_r,\phi_s)
\]

The truncation method of computing the eigenvalues consists of
choosing $n$ large enough and then finding the eigenvalues of the
matrix $M_n$ obtained from $M$ by considering only those
coefficients for which $r,s\leq n$. This is equivalent to finding
the spectrum of $M_n=P_nMP_n$, where $P_n$ is the orthogonal
projection onto $\cL_n=\lin\{\psi_1,...,\psi_n\}$. It is well
known that every point of $\Spec(M)$ is the limit of eigenvalues
of $M_n$ as $n\to\infty$. Spectral pollution arises from the fact
that the converse is not true.

As an example we consider the operator $H$ acting in $L^2(0,\pi)$
according to the formula
\[
Hf(x)=-\frac{\rmd^2 f}{\rmd x^2}
\]
and subject to Dirichlet boundary conditions. We then have
\[
\phi_n(x)=\sqrt{2/\pi}\sin(nx)
\]
and $\lam_n=n^2$ for $n=1,2,...$. We take $K$ to be the
differential operator
\[
Kf(x)=-\frac{\rmd}{\rmd x}\left\{ b(x)\frac{\rmd f}{\rmd
x}\right\}
\]
where $b(\cdot)$ is any non-negative bounded function on
$(0,\pi)$. If $D$ is the differentiation operator with domain
$W^{1,2}_0(0,\pi)$ then $H=D^\ast D$ and $K=D^\ast BD$, where $B$
is the operator of multiplication by $b(\cdot)$ acting in
$L^2(0,\pi)$. The range of $D$ is not $L^2(0,\pi)$ but the
subspace $\cM$ of all functions $f$ such that
\[
\int_0^\pi f(x)\rmd x=0.
\]
This corresponds to the fact that in the present context
(\ref{basis}) translates to
\[
\psi_n(x)=\sqrt{2/\pi}\cos(nx).
\]
Since $n=1,2,...$ these do not form a complete orthonormal set in
$L^2(0,\pi)$.

The bounded operator $M$ is given for $f\in\cM$ by
\[
Mf=bf-\langle bf,\psi_0\rangle \psi_0
\]
where $\psi_0(x)=\pi^{-1/2}$. It follows by trace class
perturbation techniques that $M$ and $B$ have the same essential
spectrum. This equals the closure of the set of essential values
of $b$, or the closure of its actual set of values if it is
piecewise continuous. However $M$ may also have a further
eigenvalue in any gap of its essential spectrum.
\begin{lemma}\label{ei} The isolated eigenvalues of $M$ are of multiplicity
$1$ and are the nonzero solutions of the equation
\[
\Xi (\lam):=\int_0^\pi\frac{b(x)}{b(x)-\lam}\rmd x =\pi.
\]
In particular if $b$ is the characteristic function of the
measurable set $E\subseteq (0,\pi)$ then $M$ has only one
isolated eigenvalue, namely $\lam=(\pi-|E|)/\pi$.
\end{lemma}

\Proof The task is to find all $f\in L^2(0,\pi)$ which satisfy
$\langle f,\psi_0\rangle =0$ and
\begin{equation}
bf-\langle bf,\psi_0\rangle \psi_0 = \lam f\label{eigenv}
\end{equation}
where $\lam$ is not in the essential range of $b$. We may
normalize $f$ by $\langle bf,\psi_0\rangle=1$ since if the LHS
equals $0$, we see from (\ref{eigenv}) that $\lam$ lies in the
essential range of $b$. We then insert
\[
f=(b-\lam)^{-1}\psi_0
\]
into the normalization equation. The function $\Xi$ is analytic
and has a positive derivative in each gap in the spectrum of $M$,
so there is at most one eigenvalue in each gap.

The entries of the matrix $A_n=P_nMP_n$ of (\ref{B}) are given for
$1\leq r,s\leq n$ by
\[
A_{n,r,s}=\frac{2}{\pi}\int_0^\pi b(x)\cos(rx)\cos(sx)\rmd x.
\]
The matrix entries of $D_n=P_nM^\ast MP_n$ are calculated using
the formula
\[
\langle M^\ast Mf,g\rangle =\int_0^\pi
b(x)^2f(x)\overline{g(x)}\rmd x-\langle bf,\psi_0\rangle\,\langle
\psi_0,bg\rangle,
\]
valid for all $f,g\in\cM$.

We consider the exactly soluble example
\[
b(x)=\left\{\begin{array}{ll}
1&\mbox{if $0\leq x\leq \alp $}\\
0&\mbox{otherwise.}
\end{array}\right.
\]
where $0 <\alp <\pi$ is arbitrary. The operator $M$ then has
two-point essential spectrum $\{0,1\}$ together with an eigenvalue
$m=(\pi -\alp)/\pi$ of multiplicity $1$. If one puts $\alp=\pi/2$
then $m=1/2$; numerical calculations show that as $n$ increases
$A_n$ has an increasing number of eigenvalues moving inwards from
$0$ and $1$ and apparently converging slowly to $0.5$. The
following table lists the four eigenvalues of $A_n$ closest to
$0.5$.
\[
\begin{array}{ccccc}
n&\lam_1&\lam_2&\lam_3&\lam_4\\
50&0.0526&0.4632&0.5368&0.9474\\
100&0.0691&0.4750&0.5250&0.9309\\
200&0.0856&0.4830&0.5170&0.9144\\
400&0.1019&0.4884&0.5116&0.8981
\end{array}
\]

The use of the techniques of this paper not only avoids the
problem of spectral pollution, but also provides accurate results
for much smaller values of $n$. We start by putting $\alp=\pi/2$
and $n=10$. The function $F_{10}(\lam)$ has three local minima.
They are given in the following table, in which the numbers are
all accurate to eight decimal places.

\[
\begin{array}{cc}
\lam&F_{10}(\lam)\\
0.00&0.00000679\\
0.50&0.10049280\\
1.00&0.00000679
\end{array}
\]

The central minimum is exactly equal to $1/2$ as a result of the
unusual symmetry for $\alp=\pi/2$. The hypothesis (H) is
therefore valid with $\nu_0=0$, $\mu_1=0.39950720$,
$\nu_1=0.60049280$ and $\mu_2=1$. An application of Lemma
\ref{induct} now yields $s_1=0.25$ and $t_1=0.75$ to six decimal
places, which confirms that the eigenvalue equals $0.5$ to the
same accuracy.

We repeat the calculation for $\alp=1.0$, in order to see whether
the unexpected accuracy of the results are caused by the
symmetry. Two representative graphs are shown
in Figure 1. 
The three local minima of $F_{8}(\lam)$ are now as follows,
accurate to six decimal places.
\[
\begin{array}{cc}
\lam&F_{8}(\lam)\\
0&0\\
0.697669&0.106699\\
0.999313&0.026051
\end{array}
\]
The hypothesis (H) is therefore valid with $\nu_0=0$,
$\mu_1=0.590970$, $\nu_1=0.804368$ and $\mu_2=1$. An application
of Lemma \ref{induct} now yields $s_1=0.340845$ and $t_1=0.840844$
to six decimal places. 
This implies that
\[
0.681688\leq \mu\leq 0.681690.
\]

The accuracy of the upper and lower bounds on the eigenvalue
$\mu$ depends on the rate of convergence in (\ref{convrate}),
which is closely associated to the rate of convergence in
\[
\lim_{n\to\infty}F_n(s)= F(s),{\hskip 0.2in}
\lim_{n\to\infty}F_n(t)=F(t)
\]
where $s=m/2$ and $t=(m+1)/2$. For reasons which we do not
understand the rate of convergence at these two points is much
faster than at other points in $(0,1)$. This phenomenon is
illustrated in the following table, in which $u$ and $v$ are
`general' points.

Putting $\alp=1.0$ we have $m=0.681690$, $s=0.340845$ and
$t=0.840845$. We also put $u=0.5$ and $v=0.75$.
\[
\begin{array}{ccccc}
n&F_n(s)&F_n(t)&F_n(u)&F_n(v)\\
2&0.340873 &0.269294& 0.254218&0.433012\\
4&0.340845 &0.162575&0.238589&0.146875\\
6&0.340845 &0.159219&0.221431&0.137537\\
8&0.340845 &0.159156&0.210113&0.116463\\
10&0.340845 &0.159155&0.207660&0.108933\\
100&0.340845 &0.159155&0.187173&0.077840\\
\infty&0.340845 &0.159155&0.181690&0.068310
\end{array}
\]
The final row lists the values of $F(\cdot)$.

The surprisingly rapid convergence of $F_n(s)$ and $F_n(t)$ does
not occur only for the two values of $\alp$ presented. A further
computation indicates that if we put $s=m/2$ then
\[
s=F(s)\leq F_{10}(s)\leq F(s)+10^{-6}
\]
for all $0<\alp<2.3$. For larger $\alp$ the difference
$F_{10}(s)-F(s)$ increases steadily. We also have
\[
s=F(s)\leq F_{50}(s)\leq F(s)+10^{-8}
\]
for all $0<\alp<2.9$.

\section{An Example with Variable Coefficients}

It might be thought that the cause of the spectral pollution in
the above examples was the fact that the basis chosen was badly
adapted to the operator.

In this section we consider the same type of example as in the
last, but using a different type of basis, constructed using the
Krylov subspace method. This example does not display spectral
pollution, but this is \emp{not} because the associated matrix is
tridiagonal. A recent example of Denisov and Simon shows that
spectral pollution can occur in the following sense even for
tridiagonal (Jacobi) matrices, \cite{densim}. They construct a
bounded, self-adjoint, tridiagonal operator $A$ acting on
$l^2(\Z^+)$ with spectrum $[-5,-1]\cup [1,5]$. The existence of a
gap in its spectrum is not a priori obvious. For certain values
of $n$, the truncation $A_n$ of $A$ to the space $\cL_n$ of
sequences with support in $\{1,n\}$ has a single eigenvalue in
$(-1,1)$. The set of limit points of these eigenvalues is equal
to the whole interval $[-1,1]$.

We replace the interval $(0,\pi)$ by $(-2,2)$ and put
\[
b(x)=\left\{ \begin{array}{ll} \alp^-+\bet^-x&\mbox{ if } x<0\\
\alp^++\bet^+x&\mbox{ if }x>0.
\end{array}\right.
\]
The essential spectrum of $b$, and of $M$, is thus
$[\alp^--2\bet^-,\alp^-]\cup[\alp^+,\alp^++2\bet^+]$. Assuming
that $\bet^\pm>0$ and $\alp^-<\alp^+$, there is an eigenvalue
$\lam$ in the spectral gap $(\alp^-,\alp^+)$. This is evaluated
by solving
\begin{equation}
\Xi (\lam):=\int_{-2}^2\frac{b(x)}{b(x)-\lam}\rmd x
=4.\label{getit}
\end{equation}
as in Lemma \ref{ei}.

We next write down the basis which we will use. Let
$\{\cP_n\}_{n=1}^\infty$ denote the Legendre polynomials, and
define
\[
Q_n(x)=\sqrt{\frac{2n+1}{2}}\cP_n(x)
\]
so that
\[
xQ_n(x)=\gam_{n+1}Q_{n+1}(x)+\gam_nQ_{n-1}(x)
\]
for all $n\geq 0$, where
\[
\gam_n=\frac{n}{\sqrt{4n^2-1}}.
\]
Then define $Q_n^\pm$ on $(-2,2)$ by
\begin{eqnarray*}
Q_n^+(x)=\left\{ \begin{array}{ll} Q_n(x-1)&\mbox{ if } x>0\\
0&\mbox{ if }x<0,
\end{array}\right. \\
Q_n^-(x)=\left\{ \begin{array}{ll} 0&\mbox{ if } x>0\\
Q_n(x+1)&\mbox{ if }x<0.
\end{array}\right.
\end{eqnarray*}
Finally define $\psi_n$ for all $n\in\Z$ by
\[
\psi_n=\left\{ \begin{array}{ll} Q_n^+&\mbox{ if } n\geq 1\\
Q_{|n|}^-&\mbox{ if }n\leq -1\\
(Q_0^++Q_0^-)/\sqrt{2}&\mbox{ if }n=0.
\end{array}\right.
\]
It is easy to check that $\{\psi_n\}_{n\in\Z}$ is a complete
orthonormal set in
\[
\cM=\{f\in L^2(-2,2):\int_{-2}^2f(x)\rmd x=0\}.
\]
Let $P_n$ be the orthogonal projection of $\cM$ onto $\cL_n=\lin\{
\psi_{-n},...,\psi_n\}$. Then it follows from standard properties
of the Legendre polynomials that $M(\cL_n)\subseteq \cL_{n+1}$.
Hence $M$ has a tridiagonal matrix. Routine calculations yield
\[
M_{r,r}=\left\{\begin{array}{ll} \alp^++\bet^+&\mbox{ if } r\geq 1\\
\alp^--\bet^-&\mbox{ if } r\leq -1\\
(\alp^++\bet^++\alp^--\bet^-)/2&\mbox{ if }r=0.
\end{array}\right.
\]

Also
\[
M_{r,r+1}=M_{r+1,r}=\left\{\begin{array}{ll}
\bet^+\gam_{r+1}&\mbox{ if } r\geq 1\\
\bet^+\gam_{1}/\sqrt{2}&\mbox{ if } r=0\\
\bet^-\gam_{1}/\sqrt{2}&\mbox{ if } r=-1\\
\bet^-\gam_{|r|}&\mbox{ if } r\leq -2.
\end{array}\right.
\]
All other coefficients of $M$ vanish. We observe that $M_{r,r}$
and $M_{r,r+1}$ have different limits as $r\to\pm\infty$.
Operators whose coefficients have different asymptotics on the
left and right have an interesting spectral and scattering theory,
studied, for example, in \cite{davsim}.

We finally compute $B_n(\lam)$ from (\ref{B}), obtaining the
formula
\[
B_n(\lam)=(A_n-\lam I)^2+E_n
\]
where $E_n=P_nM(I-P_n)MP_n$. Explicitly we have
\[
E_{n,r,s}=\left\{\begin{array}{ll}
(\bet^-\gam_{n+1})^2&\mbox{ if $r=s=-n$}\\
(\bet^+\gam_{n+1})^2&\mbox{ if $r=s=n$}\\
0&\mbox{ otherwise.}
\end{array}\right.
\]

We carried out numerical calculations for the case $\alp^+=1$,
$\alp^-=-1$, $\bet^+=40$, $\bet^-=30$, $n=50$, so that
$\dim(\cL_n)=2n+1$. The essential spectrum of $M$ is $[-61,-1]\cup
[1,81]$ and by solving (\ref{getit}) we find that $M$ also has a
single eigenvalue $\mu\sim 0.453261434$. Plotting the function
$F(\lam)$ on $[0,1]$ shows that it has a well-defined local
minimum near this point and nowhere else in the interval. From
the values
\begin{eqnarray*}
F_{50}(0.23)&=&0.22326196\\
F_{50}(0.72)&=&0.26673868
\end{eqnarray*}
we obtain the enclosures
\[
0.72-F_{50}(0.72)\leq \mu\leq 0.23+F_{50}(0.23)
\]
or, more explicitly $\mu=0.45326_{13}^{20}$. Of course all these
computations could be carried out using interval arithmetic, but
we did not do so. These upper and lower bounds on $\lam$ compare
with the completely uncontrolled values for $\lam$ obtained by the
truncation method.
\[
\begin{array}{cc}
n&  \mu\\
10&   0.50275311\\
20 &   0.45504606\\
50 & 0.45326153\\
75 &  0.45326143\\
100 &  0.45326143
\end{array}
\]
In this example we have found that by evaluating $F_{50}(\lam)$
at only two points, we obtained upper and lower bounds on $\mu$
which were of comparable accuracy to the uncontrolled value
obtained by truncation. The full plot of $F_n(\lam)$ gives
valuable geometrical insights, but may be carried out for a
smaller value of $n$ than is used for the final bounds; this is,
of course, much faster.

\vskip 0.3in
{\bf Acknowledgments} We would like to thank
M~Levitin and E~Shargorodsky for valuable discussions. \vskip
0.5in
\newpage

\newpage
\begin{figure}[h]
\begin{center}
\scalebox{0.7}{\includegraphics{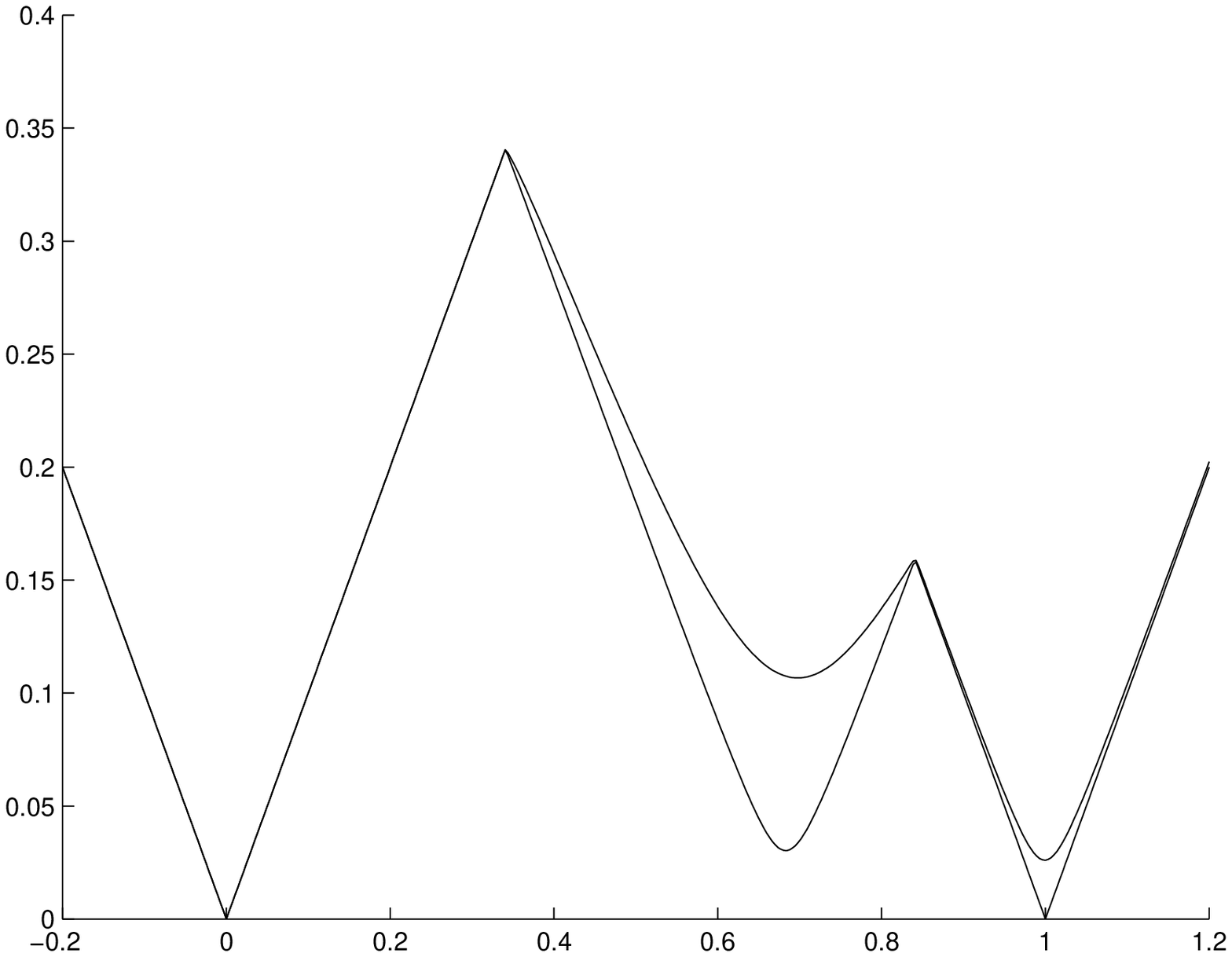}}
\end{center}
\caption{$F_{8}(\lam)$ and $F_{100}(\lam)$ for $\alp=1.0 $}
\end{figure}

\end{document}